\newcommand{\domek}{\rule{0.09in}{0.09in}}
\newcommand{\eps}{{\epsilon}_{M}}

\newenvironment{AMS}{\small\it AMS subject classification:}{ }
\newenvironment{keywords}{\small\it Keywords:}{ }
\newtheorem{theor}{Theorem}

\newtheorem{defi}{Definition}
\newtheorem{wniosek}{Proposition}

\documentclass[a4paper,12pt]{article}
\usepackage{amsfonts}
\usepackage[dvips]{graphicx}
\usepackage{latexsym}
\usepackage{ams}
\begin{document}
\author{Alicja Smoktunowicz, Przemys\l aw  Kosowski
and Iwona Wr\'{o}bel}
\title{ \Large\bf How to overcome the numerical instability of the scheme of
divided differences?}

\date{\small Faculty of Mathematics and Information Science, Warsaw University of Technology,
Pl. Politechniki 1, 00-661 Warsaw, Poland}

\maketitle

\begin{abstract}
The scheme of divided differences  is widely used in many
approximation and interpolation problems. Computing the Newton
coefficients of the interpolating polynomial is the first step of
the Bj\"{o}rck and Pereyra algorithm for solving Vandermonde
systems of equations (Cf. \cite{bjorck: 70}). Very often  this
algorithm produces very  accurate solution. The problem of
determining the Newton coefficients is intimately related with the
problem of evaluation the Lagrange interpolating polynomial, which
can be realized by many algorithms. For these reasons we use the
uniform approach and analyze also Aitken's algorithm of the
evaluation of an interpolating polynomial. We propose new
algorithms that are always numerically stable with respect to
perturbation in the function values and more accurate than the
Aitken's algorithm and the scheme of divided differences, even for
complex data.
\end{abstract}

\noindent
\begin{AMS}
{\small 65D05, 65G50}
\end{AMS}

\noindent

\section{Introduction}
The aim of this paper is to study the numerical stability of
algorithms of evaluating the Lagrange interpolating polynomial
$w_N(z)=\sum_{j=0}^{N}{f_{j} \prod_{{i=0,i\neq
j}}^{N}{\frac{z-z_{i}}{z_{j}-z_{i}}}}$ and computing the Newton
coefficients $c_0,\dots, c_N$ of $w_N(z)$. We  assume that $z_0,
z_1,\dots, z_N$ are pairwise distinct data points and $f_0,
f_1,\dots, f_N$ are associated values (real or complex) of
function $f$. Then $w_{N}(z)$ interpolates function $f(z)$ at
points $z_0, z_1, \dots, z_N$, i.e. $f(z_j)= f_j = w_N(z_j)$ for
$j=0, 1, \dots, N$ and
$w_N(z)=\sum_{j=0}^{N}{c_{j}\prod_{i=0}^{j-1}{(z-z_{i})}}$.

The problem of determining the Newton coefficients is intimately
related with the problem of evaluation the Lagrange interpolating
polynomial, which can be realized by many algorithms. For these
reasons we use the uniform approach and consider the following
general problem.

{\bf Problem I.} For given $t\in \mathbb C$, $z_0, z_1, \dots, z_N
\in \mathbb C$ (distinct) and $f_0, f_1, \dots, f_N$ compute
\begin{equation}
p_{n}(z;t)=\sum_{j=0}^{n}f_j \prod_{i=0, \hskip0.5ex i\neq j}^n
\frac{z-tz_{i}}{z_{j}-z_{i}}, \hskip3ex n=0,1,\dots,N.
\end{equation}

In the special case of $z=tz_j$ for some $j$ we obtain
$p_n(z,t)=f_j\,t^n$.

Notice that for $t=1$ we have
\begin{equation}
p_n(z;1)=w_n(z)=\sum_{j=0}^{n}f_j \prod_{i=0, \hskip0.5ex i\neq
j}^n \frac{z-z_{i}}{z_{j}-z_{i}}, \hskip3ex n=0,1,\dots,N
\end{equation}
and for $t=0$, $z=1$ we get
\begin{equation}
p_n(1;0)=c_{n}=\sum_{j=0}^{n}f_j \prod_{i=0, \hskip0.5ex i\neq
j}^n \frac{1}{z_{j}-z_{i}}, \hskip3ex n=0,1,\dots,N.
\end{equation}

\medskip
Usually the Newton coefficients are computed by the scheme of
divided differences. However, it is known (Cf. \cite{mc: 84}) that
the reputation of this algorithm is not so good. Its numerical
properties depend on the distribution and ordering of the
interpolation points. A. Kie\l basi\'nski (\cite{asia: 81}, pp.
51-53) proved that if the nodes are real and monotonically
ordered, i.e. $z_0 < z_1 < \dots < z_N$ or $z_0 > z_1 > \dots >
z_N$, then every Newton coefficient computed in fl is the exact
divided difference for slightly perturbed values $f_j$ (see
Proposition 1), that is numerically stable in the following sense.

\begin{defi} We say that an algorithm $W$ of computing
\[
p_n(z;t)=\sum_{j=0}^{n}f_j\prod_{i=0, \hskip0.5ex i\neq j}^n
\frac{z-tz_{i}}{z_{j}-z_{i}},\hskip4ex n=0,1,\dots,N
\]
is numerically stable (with respect to the data $f_0, f_1, \dots,
f_N$) if the values $\tilde p_n(z;t)$ computed by $W$ in floating
point arithmetic satisfy
\begin{equation}
\tilde p_{n}(z;t)=\sum_{j=0}^{n}[f_j(1+\delta_j^{(n)})]\prod_{i=0,
\hskip0.5ex i\neq j}^n \frac{z-tz_{i}}{z_{j}-z_{i}},\hskip5ex
\vert \delta_j^{(n)} \vert \leq \eps k_N
\end{equation}
and $k_N$ is a modest constant.
\end{defi}

It is easy to check that (4) is equivalent to
\begin{equation}
\vert \tilde p_{n}(z;t)-p_{n}(z;t)\vert \leq \eps k_N
\sum_{j=0}^{n}\vert f_j \vert \,\prod_{i=0, \hskip0.5ex i\neq j}^n
\frac{\vert z-tz_{i}\vert}{\vert z_{j}-z_{i}\vert}, \hskip3ex
n=0,\dots,N.
\end{equation}
For all set of values we have
\begin{equation}
\frac{\displaystyle{\max_{n=0,\dots,N}\vert \tilde
p_{n}(z;t)-p_{n}(z;t)\vert}}{\displaystyle{\max_{n=0,\dots,N}\vert
p_{n}(z;t) \vert}} \leq \eps k_N cond_f(z_0,\dots,z_N;z;t),
\end{equation}
where
\begin{equation}
cond_f(z_0,\dots,z_N;z;t)=\frac{\displaystyle{\max_{n=0,\dots,N}\sum_{j=0}^{n}\vert
f_j \vert \,\prod_{i=0, \hskip0.5ex i\neq j}^n \frac{\vert
z-tz_{i}\vert}{\vert
z_{j}-z_{i}\vert}}}{\displaystyle{\max_{n=0,\dots,N}\vert
p_{n}(z;t) \vert}}
\end{equation}
is the measure of sensitivity of Problem I with respect to the
data $f_0, f_1, \dots, f_N$ and will be referred to as the
condition number.

In Section 1 we review two standard methods, i.e. Aitken's
interpolation algorithm and the scheme of divided differences.
Section 3 presents a new efficient algorithm for computing (1).
Error analysis of these algorithms is given in Section 4. The new
algorithm is shown to be numerically stable in sense (4), even for
complex data, opposite to the scheme of divided differences and
Aitken's algorithm. We also show that these two algorithms are
numerically stable if the knots are real and monotonically ordered
(Cf. \cite{asia: 81}, pp. 51-53,  \cite{smok: 81}). Numerical
tests in \textsl{MATLAB} included in Section 5 confirm theoretical
results.

The new algorithm doesn't require any special ordering of
interpolation points, so we can try to reorder the knots to
achieve small relative error (see (6)). However, it seems to be a
hard problem to minimize a condition number of Problem I (see
(7)). It was observed that the Leja ordering (Cf. \cite{calvetti:
03}, \cite{lothar: 90}) may be a good choice. To reorder the
points $z_0, z_1, \dots, z_N$ by Leja, first choose $z_0$
satisfying $\vert z_0 \vert=\max_j\vert z_j \vert$, and then
determine $z_j$, $j=1,2,\dots,N-1$ such that
$\hskip1ex\prod_{k=0}^{j-1}\vert z_j-z_k \vert=\max_{i\ge
j}\prod_{k=0}^{j-1} \vert z_i-z_k \vert$.

Although in many cases helpful, sometimes it gives results
comparable to or even worse than monotone ordering, even for
equidistant points (see Section 5, Table 3).

\section{Algorithms}

The simplest way to compute Newton coefficients (3) is the well
known scheme of divided differences, which can be visualized as
follows:

\medskip

\begin{flushleft}
$c_j^{(0)}:=f_j, \hskip3ex j=0,1,\dots,N$
\end{flushleft}

\noindent $z_0 \hskip4ex c_0^{(0)}$ \hskip1ex
\raisebox{-1ex}[0pt][0pt]{$\searrow$} \hskip1ex
\raisebox{-2ex}[0pt][0pt]{$c_0^{(1)}$}
\medskip

\noindent $z_1 \hskip4ex c_1^{(0)}$ \hskip1ex
\raisebox{1.2ex}[0pt][0pt]{$\nearrow$}
\hspace{-3.1ex}\raisebox{-1.2ex}[0pt][0pt]{$\searrow$} \hskip1ex
\raisebox{-2ex}[0pt][0pt]{$c_1^{(1)}$} \hskip1ex
\raisebox{1.1ex}[0pt][0pt]{$\searrow$}
\hspace{-3.1ex}\raisebox{-1.3ex}[0pt][0pt]{$\nearrow$} \hskip1ex
$c_0^{(2)}$ \hskip7ex \raisebox{-3.8ex}[0pt][0pt]{$\searrow$}
\hspace{-3.05ex}\raisebox{-6.5ex}[0pt][0pt]{$\nearrow$}
\medskip

\noindent $z_2 \hskip4ex c_2^{(0)}$ \hskip1ex
\raisebox{1ex}[0pt][0pt]{$\nearrow$} \hskip1ex
\raisebox{-2.4ex}[0pt][0pt]{$\dots$} \hskip6ex \dots \hskip3.6ex
\dots \hskip5.9ex $c_0^{(N)}$

\vspace{2ex}

\noindent $\hspace{-0.2ex}\dots \hskip3.6ex \dots$

\vspace{2ex}

\noindent $z_N \hskip4ex c_N^{(0)}$ \hskip1ex
\raisebox{1.3ex}[0pt][0pt]{$\nearrow$} \hskip1ex
\raisebox{2.6ex}[0pt][0pt]{$c_{N-1}^{(1)}$} \hskip-0.3ex
\raisebox{3.9ex}[0pt][0pt]{$\nearrow$} \hskip1ex
\raisebox{5.2ex}[0pt][0pt]{$c_{N-2}^{(2)}$}

\vskip4ex

\noindent where $c_j^{(n)}$ are the divided differences of n{\sl
th} order. They are defined in the following way
\[
f(z_j,z_{j+1},\dots,z_{j+n})=c_j^{(n)}=\frac{c_j^{(n-1)}-c_{j+1}^{(n-1)}}{z_j-z_{j+n}}
\]
for $j=0,1,\ldots,N-n$, and $n=1,2,\ldots,N$.

Computing the value of the Lagrange interpolating polynomial at a
given point $z$ (see formula (2)) can be realized by Aitken's
algorithm (see \cite{dahlquist: 74}, Problem 8 on p. 289).

The following Algorithm I can be used for computing both Newton
coefficients and the value of interpolating polynomial, depending
on $z$ and $t$. For $t=0$ and $z=1$ algorithm reduces to the
divided differences scheme, and for $t=1$ it is Aitken's
algorithm. For all other choices of $z$ and $t$ Algorithm 1 is a
general way of evaluating (1).

\bigskip

\noindent {\large \bf Algorithm I}\\

For given (real or complex) $z_n$ and $f_n$, $t$ and $z$ this
algorithm computes
\[
p_{n}(z;t)=\sum_{j=0}^{n}f_j \prod_{i=0, \hskip0.5ex i\neq j}^n
\frac{z-tz_{i}}{z_{j}-z_{i}}, \hskip3ex n=0,1,\dots,N.
\]

\begin{flushleft}
Choose $t$ and $z$\\
$c_j^{(0)}:=f_j, \hskip3ex j=0,1,\dots,N$\\
for \hskip1ex $n=1,2, \ldots, N$\\
\hskip3.5ex for \hskip1ex $k=0,1,\ldots, N-n$
\begin{equation}
\hspace*{-19ex}c_k^{(n)}:=\displaystyle{\frac{(z-tz_{n+k})\,c_k^{(n-1)}-(z-tz_k)\,c_{k+1}^{(n-1)}}
{z_k-z_{k+n}}}
\end{equation}\\
\hskip3.5ex end\\
\hskip3.5ex $p_n(z;t):=c_{0}^{(n)}$ \\
end \\
\end{flushleft}

Then
\begin{equation}
c_k^{(n)}=\sum_{j=k}^{k+n}f_j \prod_{i=k, \hskip0.5ex i\neq
j}^{k+n} \frac{z-tz_{i}}{z_{j}-z_{i}}, \hskip3ex
n=0,1,\dots,N,\hskip3ex k=0,1, \ldots, N-n.
\end{equation}
It means that for $t=1$ $c_{k}^{(n)}$ is the Lagrange
interpolating polynomial based on the knots $z_k, z_{k+1},\dots,
z_{k+n}$ and for $z=1$, $t=0$ $c_{k}^{(n)}$ is the divided
difference of n{\sl th} order based on the knots
$z_k,\dots,z_{k+n}$.

The complexity of the divided differences and Aitken's algorithms
is $\mathcal O(N^2/2)$ and $\mathcal O(3N^2/2)$ complex
multiplications, respectively. The general algorithm requires
$\mathcal O(5N^2/2)$ complex multiplications.

\section{New algorithm}

In this section we propose a new algorithm for computing (1). The
idea is similar in spirit to the ones introduced in \cite{smok:
81}, \cite{wer: 84} and \cite{reimer: 82}, but our Algorithm II is
more efficient and robust (see Theorem 2).

For $z\neq t z_j$, $j=0,1,\dots, N$
\[
\prod_{i=0, \hskip0.5ex i\neq j }^n
(z-tz_i)=\,\,\frac{A_n}{z-tz_{j}}, \hskip3ex n=0,1,\dots,N,
\]
where
\begin{equation}
A_n=(z-tz_0)(z-tz_1)\cdot\dots\cdot(z-tz_n).
\end{equation}
Thus
\begin{equation}
p_{n}(z;t)=A_n\sum_{j=0}^{n}\frac{f_j}{(z-tz_{j})\,w_j^{(n)}},
\hskip3ex w_j^{(n)}=\prod_{i=0, \hskip0.5ex i\neq j}^n (z_j-z_i).
\end{equation}

Note that $A_n$ has the following form
\[
A_n=\left\{\begin{array}{ll} (z-z_0)\cdot\dots\cdot(z-z_n) &
\mathrm{for}\hskip1ex t=1,\\
1 & \mathrm{for} \hskip1ex t=0, \,z=1. \end{array} \right.
\]
We observe that
\[
w_j^{(0)}=1, \hskip4ex w_j^{(n)}=w_j^{(n-1)}(z_j-z_n) \hskip3ex
\mathrm{for}\hskip1ex j=0,1,\dots, n-1,
\]
\[
A_n=(z-tz_n)A_{n-1},\hskip2ex n=0,1,\dots, N,
\hskip2ex\mathrm{with}\hskip1ex A_{-1}=1.
\]
Then for $n=0,1,\dots,N$ we have $p_{n}(z;t)=A_n\,\,B_n$ where
\begin{equation}
B_n=\, \sum_{j=0}^{n}b_j^{(n)}, \hskip4ex
b_j^{(n)}=\frac{f_j}{(z-tz_{j})\,w_j^{(n)}}, \hskip3ex
j=0,1,\dots,n.
\end{equation}
We see that
\begin{equation}
b_j^{(0)}=\frac{f_j}{z-tz_{j}}, \hskip3ex j=0,1,\dots,N
\end{equation}
and
\begin{equation}
b_j^{(n)}=\frac{b_j^{(0)}}{w_j^{(n)}}, \hskip3ex j=0,1,\dots,n.
\end{equation}
Notice that if $B_{n-1}$ is written in the form
\[
B_{n-1}=b_0^{(n-1)}+b_1^{(n-1)}+\dots+b_{n-1}^{(n-1)},
\]
then we can rewrite $B_n$ as follows
\[
B_n=\frac{b_0^{(n-1)}}{(z_0-z_n)}+\frac{b_1^{(n-1)}}{(z_1-z_n)}+\dots+\frac{b_{n-1}^{(n-1)}}{(z_{n-1}-z_n)}+
\frac{b_n^{(0)}}{(z_n-z_0)(z_n-z_1)\!\dots\!(z_n-z_{n-1})}.
\]

Further, the formulae (11)-(12) let us write
\begin{equation}
b_j^{(n)}=\displaystyle{\frac{b_j^{(n-1)}}{z_j-z_n}}, \hskip3ex
j=0,1,\dots,n-1, \hskip2ex n=1,2,\dots, N
\end{equation}
and
\begin{equation}
b_n^{(n)}=\displaystyle{\frac{b_n^{(0)}}{\prod_{j=0}^{n-1}(z_n-z_j)}}.
\end{equation}

\bigskip

\noindent {\large \bf Algorithm II} \nopagebreak

For given (real or complex) $z_n$ and $f_n$, $t$ and $z\neq tz_j$
this algorithm computes
\[
p_{n}(z;t)=\sum_{j=0}^{n}f_j \prod_{i=0, \hskip0.5ex i\neq j}^n
\frac{z-tz_{i}}{z_{j}-z_{i}}, \hskip3ex n=0,1,\dots,N.
\]

\begin{flushleft}
Choose $t$ and $z$\\
$b_j^{(0)}:=\displaystyle{\frac{f_j}{z-tz_j}}, \hskip3ex j=0,1,\dots,N,$\\
$A_0:=z-tz_0$\\
$B_0:=b_0^{(0)}$\\
for \hskip1ex $n=1,2, \ldots, N$\\
\hskip3.5ex $A_n:=(z-tz_n)A_{n-1}$ \\
\hskip3.5ex for \hskip1ex $j=0,1,\ldots, n-1$\\
\hskip7ex $b_j^{(n)}:=\displaystyle{\frac{b_j^{(n-1)}}{z_j-z_n}}$ \\
\hskip3.5ex end\\
\hskip3.5ex $b_n^{(n)}:=\displaystyle{\frac{b_n^{(0)}}{\prod_{j=0}^{n-1}(z_n-z_j)}}$ \\
\hskip3.5ex $B_n:=\sum_{j=0}^n b_j^{(n)}$ \\
\hskip3.5ex $p_n(z;t):=A_n B_n$ \\
end \\
\end{flushleft}

The complexity of Algorithm II is $\mathcal O(N^2)$ flops.

Note that this algorithm can be implemented more efficiently than
it's written now, namely $b_n^{(n)}$ can be computed along with
$b_j^{(n)}$ in the inner loop for
\[
b_n^{(n)}:=\frac{(-1)^n\,b_n^{(0)}}{\prod_{j=0}^{n-1}(z_j-z_n)}.
\]
We decided to write Algorithm II this way to make it more clear.

\vskip5ex

\section{Error analysis of algorithms}

We consider complex arithmetic (cfl) (Cf. \cite{smok: 81})
implemented using real arithmetic (fl) with machine unit $\eps$.
We assume that for $x, y \in \mathbb R$ we have
\begin{equation}
fl(x \pm y) = (x \pm y) (1+\delta), \hskip3ex \vert \delta \vert
\leq \eps.
\end{equation}
Then
\begin{equation}
cfl(x \odot y) = (x \odot y) (1+\delta), \hskip3ex \vert \delta
\vert \leq c_\odot\,\eps \hskip3ex \mathrm{for} \;\;x,y \in
\mathbb C,
\end{equation}
where
\begin{equation}
c_\odot=\left\lbrace\begin{array}{cl} 1 & \mathrm{for} \hskip.5cm
\odot \in \{+,-\},\\
1+\sqrt 2 & \mathrm{for} \hskip.5cm \odot=*,\\
4+\sqrt 2 & \mathrm{for} \hskip.5cm \odot=/.
\end{array}\right.
\end{equation}
The constant $c_\odot$ is not significant, it depends on the
implementation of floating point operations. In fact, the bounds
of errors in Algorithms I and II are of the same form for real and
complex arithmetic, only the constants are increased
appropriately.

Our error analysis is uniform, it includes both real and complex
cases, i.e. the data $z_j$, $f_j$, $z$ and $t$ can be either real
or complex.

The values computed in fl or cfl will be marked with tilde in
order to distinguish them from the exact ones.

\subsection{Error analysis of Algorithm I}

\begin{theor} Assume that $z_j$, $f_j, \in \mathbb C$ for $j=0,1,\ldots, N$ and $z$, $t \in \mathbb C$ are
exactly representable in cfl and
\begin{equation}
cfl(tz_j)=tz_j, \hskip3ex j=0,1,\dots,N.
\end{equation}
Then the  values $\tilde c_{k}^{(n)}$ for $n=0,1,\dots N$,
$k=0,1,\dots,N-n$ computed in cfl by Algorithm I satisfy
\begin{equation}
\tilde
c_{k}^{(n)}=\sum_{j=k}^{k+n}[f_j(1+\delta_{k,j}^{(n)})]\prod_{i=k,
\hskip0.5ex i\neq j}^{k+n} \frac{z-tz_{i}}{z_{j}-z_{i}},\hskip3ex
\vert \delta_{k,j}^{(n)} \vert \leq \eps n \,d\, L^{n-1} +\mathcal
O(\eps^{\,2}),
\end{equation}
where
\begin{equation}
d=\left\lbrace\begin{array}{cl} 5 & in \hskip1ex the \hskip1ex real \hskip1ex case,\\
8+2 \sqrt 2 & in \hskip1ex the \hskip1ex complex \hskip1ex case
\end{array}\right.
\end{equation}
and
\begin{equation}
L=\max_{0 \leq i < j < k \leq N}\frac{\vert z_{i}-z_{j} \vert+
\vert z_{j}-z_{k} \vert}{\vert z_{i}-z_{k} \vert}.
\end{equation}
\end{theor}

\noindent {\it Proof.} By induction on $n$. It follows from (8)
and (17)-(19) that for $n=0,1,\dots N$, $k=0,1,\dots,N-n$ the
computed quantities satisfy
\begin{equation}
\tilde
c_k^{(n)}:=\displaystyle{\frac{(1+\alpha_{k}^{(n)})(z-tz_{n+k})\tilde
c_k^{(n-1)}-(1+\beta_{k}^{(n)})(z-tz_k)\tilde
c_{k+1}^{(n-1)}}{z_k-z_{k+n}}},
\end{equation}
where
\begin{equation}
\vert \alpha_{k}^{(n)} \vert, \, \vert \beta_{k}^{(n)} \vert \leq
\eps\, d +\mathcal O(\eps^{\,2}).
\end{equation}
It is clear that for $n=1$ the result is immediate, since
\[
\tilde
c_k^{(1)}=f_k\,(1+\alpha_{k}^{(1)})\frac{z-tz_{k+1}}{z_k-z_{k+1}}+
f_{k+1}\,(1+\beta_k^{(1)})\frac{z-tz_k}{z_{k+1}-z_k},
\]
so we can take $\delta_{k,k}^{(1)}=\alpha_k^{(1)}$ and
$\delta_{k,k+1}^{(1)}=\beta_k^{(1)}$, $k=0,1,\dots, N-1$.

Now assume that (21) is true for $n-1$. We will prove that it
holds also for $n$. By assumption for $n-1$ we get
\begin{equation}
(z-tz_{n+k})\,\tilde
c_{k}^{(n-1)}=\sum_{j=k}^{k+n-1}[f_j(1+\delta_{k,j}^{(n-1)})]\prod_{i=k,
\hskip0.5ex i\neq j}^{k+n} \frac{z-tz_i}{z_j-z_i}\,(z_j-z_{k+n}),
\end{equation}
\begin{equation}
(z-tz_k)\,\tilde
c_{k+1}^{(n-1)}=\sum_{j=k+1}^{k+n}[f_j(1+\delta_{k+1,j}^{(n-1)})]\prod_{i=k,
\hskip0.5ex i\neq j}^{k+n} \frac{z-tz_i}{z_j-z_i}\,(z_j-z_k).
\end{equation}
To simplify the notation let us define
\[
1+\phi_{k,j}^{(n)}=(1+\alpha_k^{(n)}) (1+\delta_{k,j}^{(n-1)})
\]
and
\[
1+\psi_{k,j}^{(n)}=(1+\beta_k^{(n)}) (1+\delta_{k+1,j}^{(n-1)}),
\]
bounded as follows
\begin{equation}
\vert \phi_{k,j}^{(n)}\vert,\, \vert \psi_{k,j}^{(n)} \vert \leq
\eps d\,(1 + (n-1) L^{n-2}) +\mathcal O(\eps^{\,2}).
\end{equation}
Subtraction of (26) from (27) leads to
\[
\tilde c_k^{(n)}=\sum_{j=k}^{k+n} f_j (1+\delta_{k,j}^{(n)})
\prod_{i=k, \hskip0.5ex i\neq j}^{k+n} \frac{z-tz_i}{z_j-z_i},
\]
where
\[
1+\delta_{k,j}^{(n)} = \frac{(1+\phi_{k,j}^{(n)})(z_j-z_{k+n})+
(1+\psi_{k,j}^{(n)})(z_k-z_j)}{z_k-z_{k+n}},
\]
for $j=k+1,\dots,k+n-1$, and
\[
1+\delta_{k,k}^{(n)} = 1+\phi_{k,k}^{(n)}, \hskip5ex
1+\delta_{k,k+n}^{(n)} = 1+\psi_{k,k+n}^{(n)}.
\]
This leads to the estimation
\[
\vert \delta_{k,j}^{(n)} \vert \leq \max\{\vert
\phi_{k,j}^{(n)}\vert, \, \vert \psi_{k,j}^{(n)} \vert\} \, L,
\hskip4ex \mathrm{for} \hskip2.5ex j=k,\dots,k+n,
\]
which, in combination with (28) and the fact that $L \ge 1$ (see
(23)) results in $\vert \delta_{k,j}^{(n)} \vert \leq \eps n
d\,L^{n-1} + \mathcal O(\eps^{\,2})$, which completes the proof.
\domek

\medskip

Although, in general, the constant $L$ defined in (23) might be
large, in the special cases the Algorithm I is numerically stable
in sense (4). This happens, for example, for specially ordered
knots. For details see the following proposition and remarks below
it.

\begin{wniosek} Assume that $z_j$, $f_j \in \mathbb R$ for $j=0,1,\ldots, N$ and $z, t\in \mathbb R$ are
exactly representable in fl and {\rm(20)} holds. If the knots
$z_j$ are monotonically ordered then the values $\tilde
p_{n}(z;t)= \tilde c_0^{(n)}$ computed in fl by Algorithm I
satisfy
\begin{equation}
\tilde
p_{n}(z;t)=\sum_{j=0}^{n}[f_j(1+\delta_j^{(n)})]\,\prod_{i=0,
\hskip0.5ex i\neq j}^n \frac{z-tz_{i}}{z_{j}-z_{i}},\hskip4ex
\vert \delta_j^{(n)} \vert \leq \eps 5\,n +\mathcal O(\eps^{\,2}).
\end{equation}
\end{wniosek}

\noindent {\it Proof.} It is a direct consequence of Theorem 1 for
$k=0$. For monotonically ordered knots $z_j$ the constant $L$ in
(23) is 1. \domek

\medskip
{\bf Remark.} If $z_j \in \mathbb C$, $z_j=a+jh$, $j=0,1,\dots,N$,
$h=(b-a)/N$, $a, b \in \mathbb C$, then $L=1$ and Algorithm I is
also numerically stable, i.e. the condition (4) holds with the
constant $k_N \approx (8+2 \sqrt 2)N$.

\medskip
In Section 5 we present numerical experiments showing that
Algorithm I is not always numerically stable.

\subsection{Error analysis of Algorithm II}

\begin{theor} Assume that $z_j$, $f_j \in \mathbb C$, $z\neq tz_j$ for $j=0,1,\ldots, N$ and $z, t\in \mathbb C$ are
exactly representable in cfl and {\rm(20)} holds. Then the values
$\tilde p_{n}(z;t)$ computed in cfl by Algorithm II satisfy
\begin{equation}
\tilde
p_{n}(z;t)=\sum_{j=0}^{n}[f_j(1+\delta_j^{(n)})]\,\prod_{i=0,
\hskip0.5ex i\neq j}^n \frac{z-tz_{i}}{z_{j}-z_{i}},\hskip5ex
\vert \delta_j^{(n)} \vert \leq \eps k_N+\mathcal O(\eps^{\,2}),
\end{equation}
where
\begin{equation}
k_N = (N+1)\left\lbrace\begin{array}{cl} 5 & in \hskip1ex the \hskip1ex real \hskip1ex case,\\
8+2\sqrt 2 & in \hskip1ex the \hskip1ex complex \hskip1ex case.
\end{array}\right.
\end{equation}
\end{theor}

\noindent {\it Proof.} The proof is straightforward. Using
(13)-(16) and (17)-(19) we obtain the following quantities,
computed in fl or cfl
\[
\tilde b_j^{(0)}=\frac{f_j}{z-tz_j}(1+\Delta_j^{(0)}), \hskip3ex
j=0,1,\ldots, N,
\]
and for $n=1,2,\ldots, N$
\[
\tilde A_n=A_n(1+\alpha_n),
\]
\[
\tilde b_j^{(n)}=\frac{\tilde
b_j^{(n-1)}}{z_j-z_n}(1+\beta_j^{(n)}), \hskip3ex j=0,1,\ldots,
n-1,
\]
\[
\tilde b_n^{(n)}=\frac{\tilde
b_n^{(0)}}{\prod_{j=0}^{n-1}(z_n-z_j)}(1+\gamma_n),
\]
From (12) we obtain
\[
\tilde B_n=\sum_{j=0}^n \tilde b_j^{(n)}(1+\alpha_j^{(n)}),
\hskip4ex \vert \alpha_j^{(n)} \vert \leq \eps\, n +\mathcal
O(\eps^{\,2}),
\]
\[
\tilde p_n(z;t)=\tilde A_n \tilde B_n (1+\delta_n),
\]
where
\[
\vert \Delta_j^{(0)} \vert,\; \vert \beta_j^{(n)} \vert \leq \eps
\left\lbrace\begin{array}{cl} 2 &
\mathrm{(real} \hskip1ex \mathrm{case)}\\
5+\sqrt 2 & \mathrm{(complex} \hskip1ex \mathrm{case)}
\end{array}\right.+\mathcal O(\eps^{\,2}),
\]
\[
\vert \alpha_n \vert \leq \eps \left\lbrace\begin{array}{cl}
2\,n+1 & \mathrm{(real} \hskip1ex \mathrm{case)}\\
(2+\sqrt 2)\,n + 1 & \mathrm{(complex} \hskip1ex \mathrm{case)}
\end{array}\right.+\mathcal O(\eps^{\,2}),
\]
\[
\vert \gamma_n \vert \leq \eps \left\lbrace\begin{array}{cl} 2\,n
&
\mathrm{(real} \hskip1ex \mathrm{case)}\\
(2+\sqrt 2)\,n + 3 & \mathrm{(complex} \hskip1ex \mathrm{case)}
\end{array}\right.+\mathcal O(\eps^{\,2}),
\]
\[
\vert \delta_n \vert \leq \eps \left\lbrace\begin{array}{cl} 1 &
\mathrm{(real} \hskip1ex \mathrm{case)}\\
1+\sqrt 2 & \mathrm{(complex} \hskip1ex \mathrm{case)}
\end{array}\right.+\mathcal O(\eps^{\,2}).
\]
Now the formulae (15), (16) yield for all $b_j^{(n)}$
\begin{equation}
\tilde b_j^{(n)}=b_j^{(n)} (1+\phi_j^{(n)}),
\end{equation}
where
\begin{equation}
\vert \phi_j^{(n)} \vert \leq \eps (n+1)
\left\lbrace\begin{array}{cl} 2 &
\mathrm{(real} \hskip1ex \mathrm{case)}\\
5+\sqrt 2 & \mathrm{(complex} \hskip1ex \mathrm{case)}
\end{array}\right.+\mathcal O(\eps^{\,2})
\end{equation}
and for $\tilde B_n$, $n=1,2\dots,N$
\[
\tilde B_n=\sum_{j=0}^n \tilde b_j^{(n)}(1+\psi_j^{(n)}),
\]
where
\[
\vert \psi_j^{(n)} \vert \leq \eps (n+1)
\left\lbrace\begin{array}{cl} 3 &
\mathrm{(real} \hskip1ex \mathrm{case)}\\
6+\sqrt 2 & \mathrm{(complex} \hskip1ex \mathrm{case)}
\end{array}\right.+\mathcal O(\eps^{\,2}).
\]
Finally, this and (32), (33) lead to (30), which proves the
theorem. \domek

\bigskip
{\bf Remark.} In the case of divided differences algorithm
(formula (3)) the constants are smaller, namely (Cf. \cite{smok:
81})
\[
\vert \delta_j^{(n)} \vert \leq \eps (n+1)
\left\lbrace\begin{array}{cl} 3 &
\mathrm{(real} \hskip1ex \mathrm{case)}\\
6+\sqrt 2 & \mathrm{(complex} \hskip1ex \mathrm{case)}
\end{array}\right.+\mathcal O(\eps^{\,2}).
\]

\bigskip

\section{Numerical experiments}

To illustrate our results we present numerical experiments carried
out in \textsl{MATLAB} with unit roundoff $\eps=2.2e\!\!-\!\!16$.

We implemented all methods and performed many tests in order to
investigate the behaviour of these methods, i.e. to compare the
errors generated by each of them. This paragraph contains the
results of some of these tests.

Algorithm II is general, two special cases are computing the
Newton coefficients and computing the value of an interpolating
polynomial. These are the cases we considered in our tests.

The first kind of tests involves computing the Newton
coefficients. One type of the functions we interpolated was
$f(z)=z^s$, $s=2,3,\dots,7$.

To compare our algorithm with the classical scheme of divided
differences we evaluated the following quantities:
\begin{equation}
error1=\frac{\displaystyle{\max_{n=s+2,\dots,N}\vert \tilde c_n
\vert}}{\displaystyle{\eps \max_{n=0,1,\dots,N}\sum_{j=0}^{n}\vert
f_j \vert \prod_{i=0, \hskip0.5ex i\neq j}^{n} \frac{1}{\vert
z_{j}-z_{i}\vert}}}
\end{equation}
and
\begin{equation}
error2=\frac{\displaystyle{\max_{n=s+2,\dots,N}\vert \tilde c_n
\vert}}{\displaystyle{\eps \max_{n=0,1,\dots,N}\vert c_n \vert}},
\end{equation}
where $\tilde c_n$ and $c_n$ are the computed and the exact Newton
coefficients, respectively. Notice that if $f(z)=z^s$, then
$c_{s+2}=\dots=c_N=0$ for $N\ge s+2$.

For numerically stable algorithms $error1$ is small, it is a
measure of instability of an algorithm, which comes from the
following chain of inequalities (see (4)-(7))
\begin{equation}
\hspace*{-3ex}\max_{n=s+2,\dots,N}\!\vert \tilde c_n \vert \! \leq
\!\max_{n=0,\dots,N}\!\vert \tilde c_n\! - c_n \vert \! \leq \!
\eps k_N \max_{n=0,\dots,N}\! \sum_{j=0}^{n}\vert f_j \vert \!
\prod_{i=0, \hskip0.5ex i\neq j}^{n} \! \frac{1}{\vert
z_{j}-z_{i}\vert}.
\end{equation}
Now it is easy to see that $error1$ should be less than or equal
to the constant $k_N$ in (4).

Tables 1 and 2 contain the respective values of $error1$ and
$error2$ for Aitken's algorithm and Algorithm II applied to the
function $f(z)=z^7$. The interpolation was based on 50 or 80
random knots ordered either monotonically or by Leja ordering. In
all three tests problem was well-conditioned, with condition
number (7) of order 1.

\begin{center}
\noindent {\footnotesize Table 1: The error (34) for the scheme of divided differences (Div. diff.) and Algorithm II} \\
\hspace*{-16.5ex} {\footnotesize for random knots,
increasing ($\uparrow$) and Leja orderings.}\\
\nopagebreak
\bigskip
\nopagebreak
\begin{tabular}{|c|c|c|c|c|}
\hline $N$ & Div. diff. + $\uparrow$ & Alg. II + $\uparrow$ & Div.
diff. + Leja & Alg. II + Leja \\
\hline $50$ & $1.161255e\!+\!04$ & $5.470269e\!-\!02$ & $4.587882e\!-\!01$ & $6.130860e\!-\!05$\\
\hline $50$ & $8.622460e\!+\!02$ & $4.401986e\!-\!02$ & $8.014843e\!-\!05$ & $3.478342e\!-\!05$\\
\hline $80$ & $2.274724e\!+\!06$ & $1.340639e\!-\!01$ & $1.195999e\!-\!04$ & $7.967939e\!-\!05$\\
\hline
\end{tabular}
\end{center}
\medskip
\begin{center}
\noindent {\footnotesize Table 2: The error (35) for the scheme of
divided differences (Div. diff.)} \nopagebreak
\\ \hspace*{-5ex}
{\footnotesize  and Algorithm II for the same data as in Table
1.}\nopagebreak
\\
\nopagebreak
\bigskip
\nopagebreak
\begin{tabular}{|c|c|c|c|c|}
\hline $N$ & Div. diff. + $\uparrow$ & Alg. II + $\uparrow$ & Div.
diff. + Leja & Alg. II + Leja \\
\hline $50$ & $2.019597e\!+\!04$ & $9.513618e\!-\!02$ & $4.587882e\!-\!01$ & $6.130860e\!-\!05$\\
\hline $50$ & $9.928299e\!+\!02$ & $5.068650e\!-\!02$ & $8.014843e\!-\!05$ & $3.478342e\!-\!05$\\
\hline $80$ & $2.278485e\!+\!06$ & $1.342856e\!-\!01$ & $1.195999e\!-\!04$ & $7.967939e\!-\!05$\\
\hline
\end{tabular}
\end{center}
\bigskip
\medskip

Now let us turn our attention to evaluation of an interpolating
polynomial.

Tests presented below were performed for the function $f(z)=z^7$.
Interpolation was based on 100 equally spaced knots from the
interval $[-1,1]$, ordered either monotonically or by Leja
ordering.

All graphs present the error (see (4)-(7) for the function
$f(z)=z^7$)
\begin{equation}
error3=\frac{\vert \tilde p_n(z,1) - z^7 \vert}{\displaystyle{\eps
\max_{n=0,1,\dots,N}\sum_{j=1}^{n}\vert f_j\vert\prod_{i=1,
\hskip0.5ex i\neq j}^{n} \frac{\vert z-z_i \vert}{\vert
z_{j}-z_{i}\vert}}}
\end{equation}
computed for equally spaced points from interval $[-98/99,97/98]$,
where $\tilde p_n(z,1)$ is the value of an interpolating
polynomial $p_n(z,1)$ in (2) given either by Aitken's algorithm or
by Algorithm II.

Figures 1 and 2 describe the results for Aitken's algorithm and
Algorithm II, respectively, for monotonically ordered knots. They
confirm theoretical results, that Aitken's algorithm is
numerically stable for monotonically ordered points (increasingly
or decreasingly), see Theorem 1.

\bigskip

\begin{center}
\includegraphics[width=10cm,height=7cm]{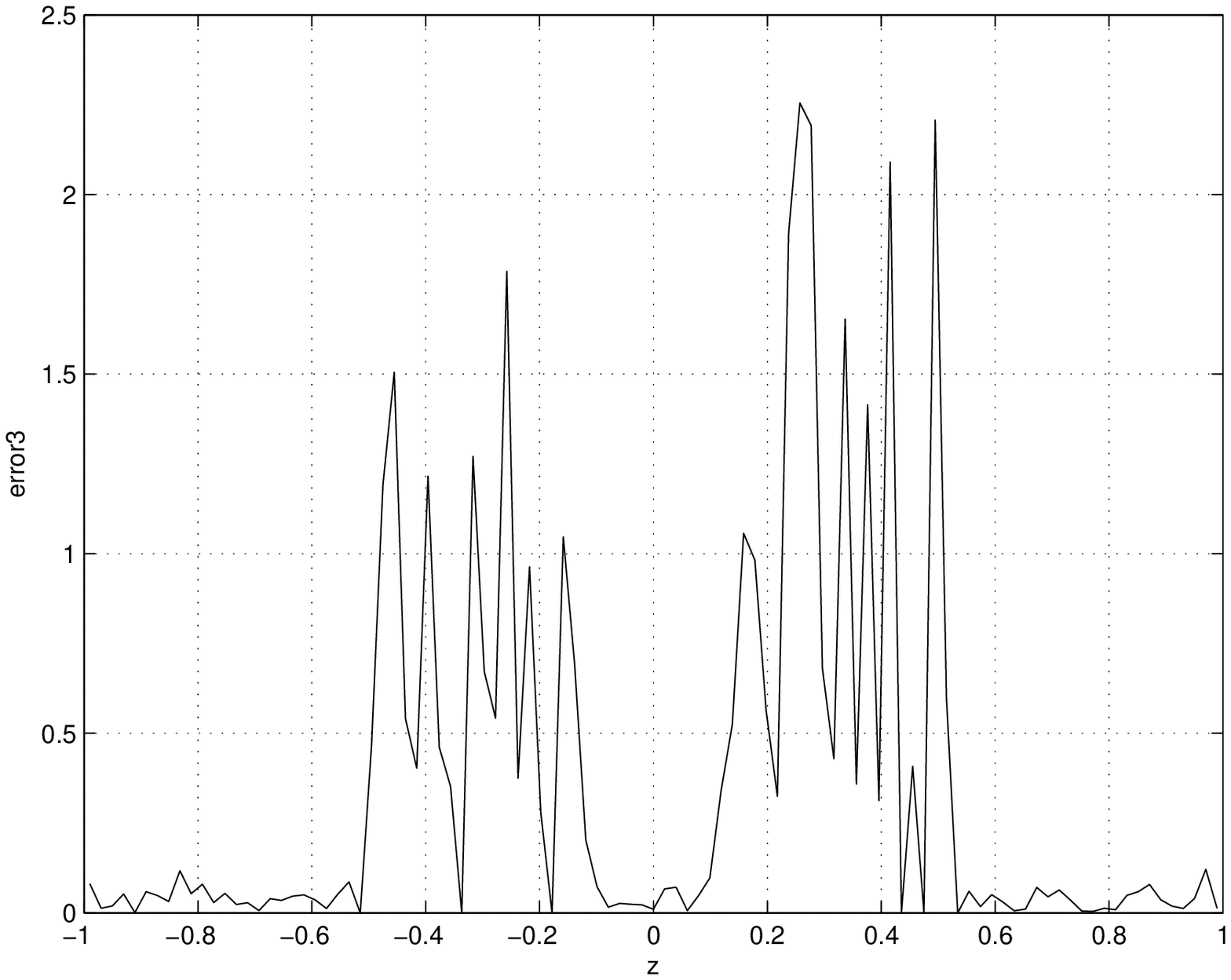}
\nopagebreak

\hspace*{-4ex}{\footnotesize Figure 1: The error (37) for Aitken's algorithm for $f(z)=z^7$}\\
\hspace*{1ex}{\footnotesize with the knots being 100 equally spaced points}\\
\hspace*{7.5ex} {\footnotesize from the interval $[-98/99,97/98]$,
increasingly ordered.}
\end{center}

\noindent
\begin{center}
\includegraphics[width=10cm,height=7cm]{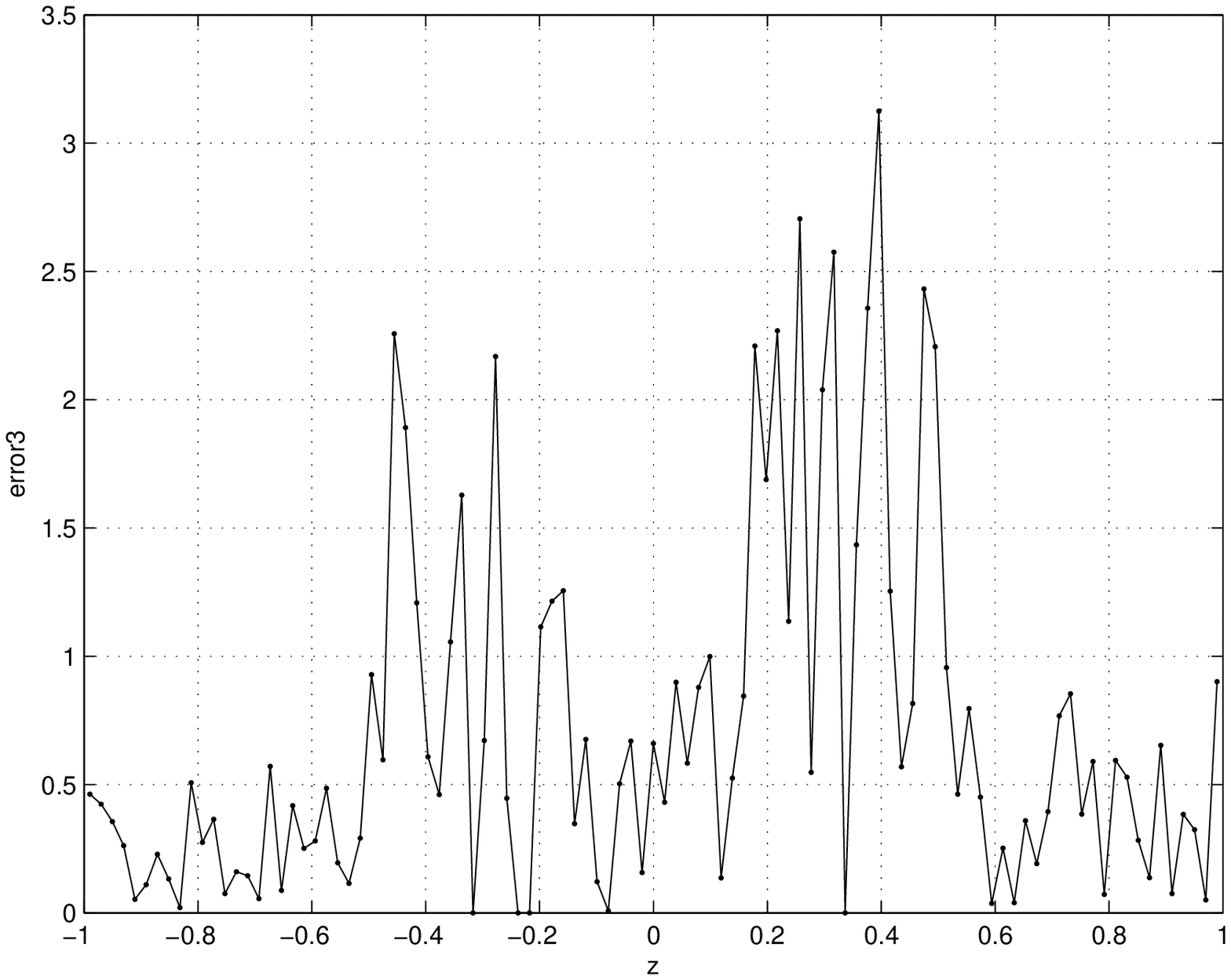}
\nopagebreak

\hspace*{-10ex}{\footnotesize Figure 2: The error (37) for Algorithm II for $f(z)=z^7$} \\
\hspace*{-0.5ex} {\footnotesize with the knots being 100 equally
spaced points}\\
\hspace*{7ex}{\footnotesize from the interval $[-98/99,97/98]$,
increasingly ordered.}
\end{center}

\bigskip Figures 3 and 4 present the values of error (37) for
Aitken's algorithm and Algorithm II, respectively, for knots
ordered by Leja ordering.

\bigskip
\begin{center}
\includegraphics[width=10cm,height=7cm]{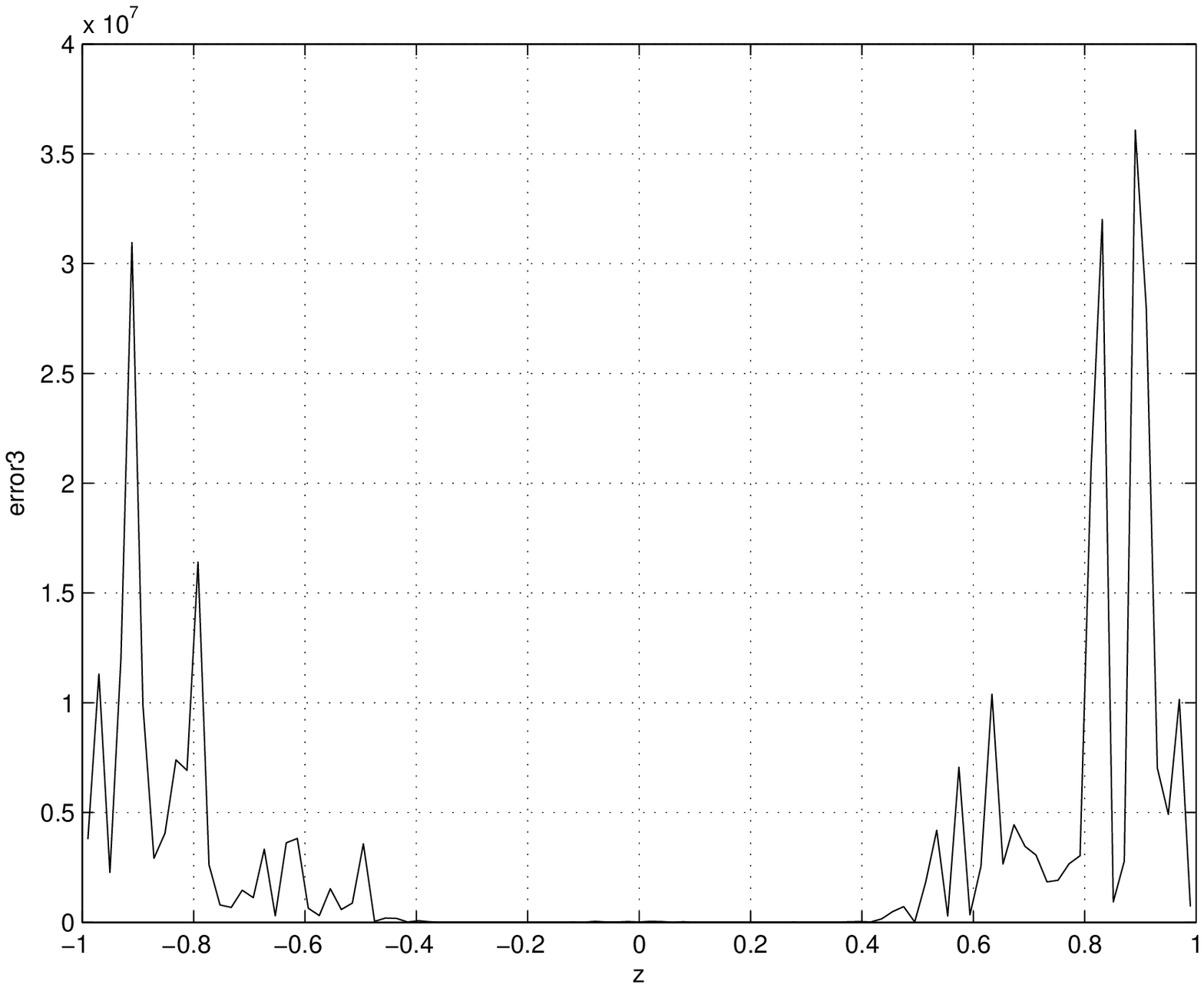}
\nopagebreak

\hspace*{-5ex}{\footnotesize Figure 3: The error (37) for Aitken's algorithm for $f(z)=z^7$}\\
{\footnotesize with the knots being 100 equally spaced points}\\
\hspace*{11ex}{\footnotesize from the interval $[-98/99,97/98]$,
ordered by Leja ordering.}
\end{center}

\noindent
\begin{center}
\includegraphics[width=10cm,height=7cm]{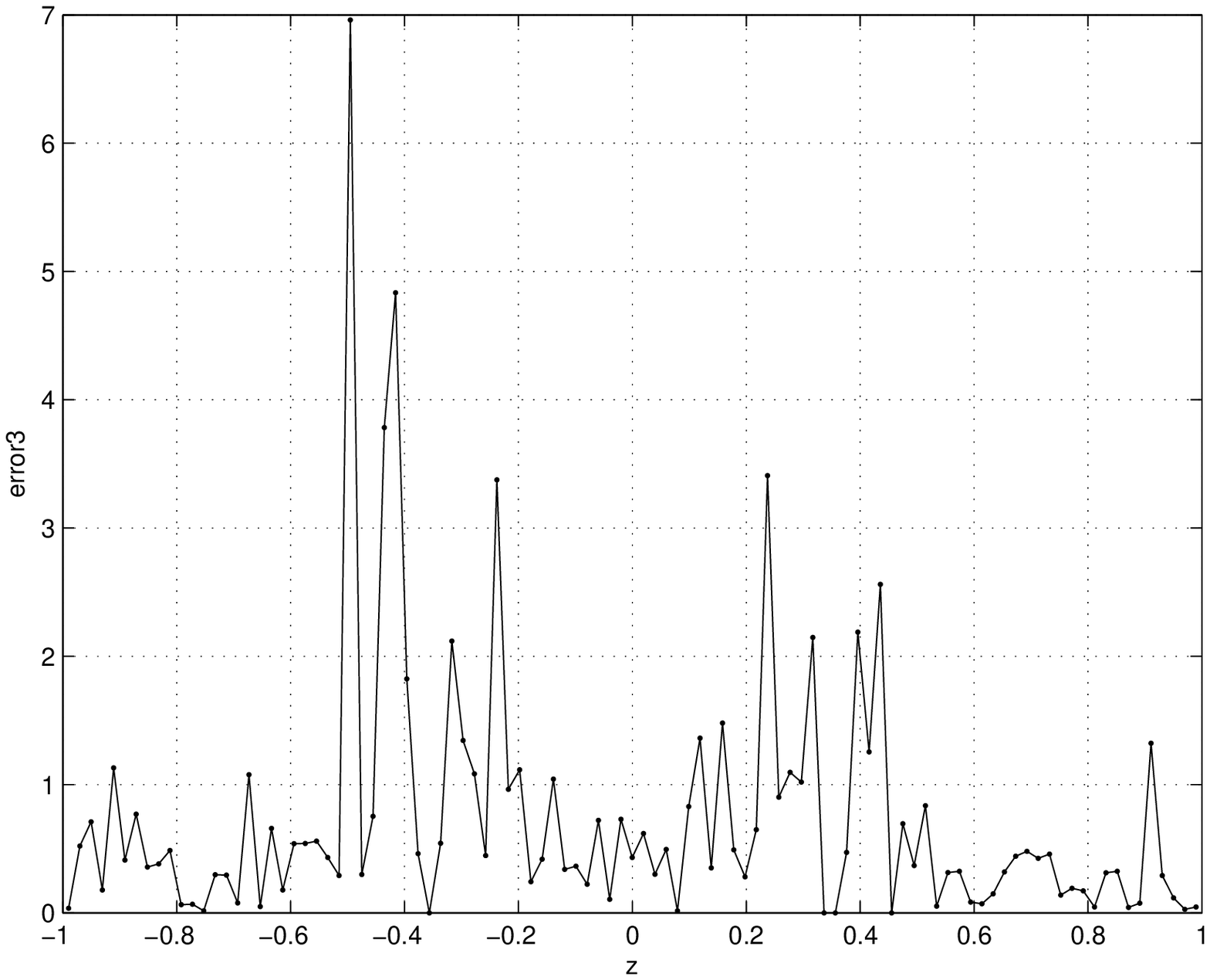}

\hspace*{-11ex} {\footnotesize Figure 4: The error (37) for
Algorithm II for $f(z)=z^7$}\\
\hspace*{-1ex}
{\footnotesize with the knots being 100 equally spaced points}\\
\hspace*{10.5ex}{\footnotesize from the interval $[-98/99,97/98]$,
ordered by Leja ordering.}
\end{center}

\bigskip

Table 3 contains values of $error3$ at five of the points in the
data used to create figures 1-4.

\medskip
\begin{center}
\noindent \hspace*{-10ex}{\footnotesize Table 3: The error (37) for Aitken's algorithm and Algorithm II} \\
\hspace*{-0.5ex} {\footnotesize with knots ordered either
increasingly or by Leja ordering.}\\
\bigskip
\begin{tabular}{|c|c|c|c|c|}
\hline $z$ & Aitken + $\uparrow$ & Alg. II + $\uparrow$ & Aitken +
Leja & Alg. II + Leja \\
\hline $-.989899$ & $8.078746e\!-\!02$ & $4.623302e\!-\!01$ & $3.791592e\!+\!06$ & $3.593606e\!-\!02$\\
\hline $-.791929$ & $7.973497e\!-\!02$ & $2.746146e\!-\!01$ & $1.640714e\!+\!07$ & $6.377801e\!-\!02$\\
\hline $-.593960$ & $3.579114e\!-\!02$ & $2.806026e\!-\!01$ & $6.421192e\!+\!05$ & $5.397304e\!-\!01$\\
\hline $-.395990$ & $1.215628e\!+\!00$ & $6.078142e\!-\!01$ & $7.446058e\!+\!04$ & $1.823442e\!+\!00$\\
\hline $-.198021$ & $2.787772e\!-\!01$ & $1.115109e\!+\!00$ & $3.299607e\!+\!03$ & $1.115109e\!+\!00$\\
\hline
\end{tabular}
\end{center}

\bigskip
Note that Leja ordering, although in many cases helpful, is not a
general remedy. Aitken's algorithm is more accurate 
for monotonically ordered knots than for Leja points. This is
because of the constant $L$ in Theorem 1. For monotone order $L$
is always 1, for Leja ordering $L$ may be big, for example in the
test presented above $L=197$.

\baselineskip=0.9\normalbaselineskip

\end{document}